\documentclass[runningheads]{llncs}
\usepackage{amsfonts, amsmath, amssymb, mathtools, bbm}
\usepackage{fancyvrb, tikzcd, extpfeil}
\usepackage{fancyhdr}

\usepackage{fancyhdr}
\renewcommand\chaptermark[1]{\markboth{\thechapter\;\;\; #1}{}}
\begingroup
    \makeatletter
    \@for\theoremstyle:=definition,remark,plain\do{%
        \expandafter\g@addto@macro\csname th@\theoremstyle\endcsname{%
            \addtolength\thm@preskip\parskip
            }%
        }
\endgroup 

\newcommand{\ra}{\rightarrow}

\newcommand{\mcc}{\mathcal{C}}

\newcommand{\mcj}{\mathcal J}

\newcommand{\mco}{\mathcal{O}}

\newcommand{\mct}{\mathcal{T}}

\newcommand{\mbo}{\mathbb{O}}

\newcommand{\vp}{\varphi} 

\newcommand{\stl}{\left\{}
\newcommand{\str}{\right\}}

\newcommand{\hm}{\text{Hom}}

\newcommand{\und}{\underline}
\newcommand{\mbbo}{\mathbbm{1}}
\newcommand{\sets}{\textbf{Sets}}

\usepackage[T1]{fontenc}
\usepackage{graphicx}
\begin{document}
\title{A Formal Algebraic Framework for DSL Composition\thanks{Supported by DARPA V-SPELLS HR001120S0058.\\ \vspace{-.25 cm}\\
\textbf{Distribution Statement A}:  Approved for Public Release, Distribution Unlimited}}
%
%
\author{Zachary Flores\and
Angelo Taranto\and
Eric Bond}
\authorrunning{Z. Flores et al.}
%
\institute{Two Six Technologies, Arlington, VA 22203, USA \\
\email{zachary.flores@twosixtech.com}\\ \email{angelo.taranto@twosixtech.com}\\ \email{eric.bond@twosixtech.com}}
\maketitle              
\begin{abstract}
We discuss a formal framework for using algebraic structures to model a meta-language that can write, compose, and provide interoperability between abstractions of DSLs.  The purpose of this formal framework is to provide a verification of compositional properties of the meta-language.  Throughout our paper we discuss the construction of this formal framework, as well its relation to our team's work on the DARPA V-SPELLS program via the pipeline we have developed for completing our verification tasking on V-SPELLS.  We aim to give a broad overview of this verification pipeline in our paper.  The pipeline can be split into four main components:  the first is providing a formal model of the meta-language in Coq; the second is to give a specification in Coq of our chosen algebraic structures; third, we need to implement specific instances of our algebraic structures in Coq, as well as give a proof in Coq that this implementation is an algebraic structure according to our specification in the second step; and lastly, we need to give a proof in Coq that the formal model for the meta-language in the first step is an instance of the implementation in the third step.   

\vspace{.25 cm}
\textbf{Disclaimer:}  The views, opinions and/or findings expressed are those of the author and should not be interpreted as representing the official views or policies of the Department of Defense or the U.S. government.
\end{abstract}

\section{Introduction}\label{sec1}
 \lfoot{\textbf{Distribution Statement A}:  Approved for Public Release, Distribution Unlimited}

The DARPA V-SPELLS (Verified Security and Performance Enhancement of Large Legacy Software) program aims to create developer-accessible capability for piece-by-piece enhancement of software components for large legacy codebases with new verified code that is \textit{safely composable }with the rest of the system.  

In our approach with the Johns Hopkins Applied Physics Laboratory to solving the problems posed by V-SPELLS, our tool in development, called LUMOS, begins by applying tools from static analysis, natural language processing, and dynamic analysis to the legacy source code in order to generate domain-specific semantic models (DSSMs) from the DSLs that comprise the source code.  These DSSMs will be generated in a language we refer to as the \textit{meta-DSL}, and in order to provide the patches to the legacy code requested in V-SPELLS, these models for the DSLs will have to be composed in very specific ways.  In order to ensure correctness of composition, we are providing verification via an algebraic framework using several ideas from category theory in which the key structure to our modeling is called an \textit{operad}.  Operads have begun to play an increasing role within applied mathematics (see \cite{LBPF}, \cite{FBSD}, \cite{BF}, \cite{BO}, and \cite{GLMN}), and we find they provide an excellent mathematical model for our verification needs on V-SPELLS.  

Our focus in this paper is to give a brief overview of both our formal modeling for DSLs, and the mathematical framework we will be using in this formalization.  With that, we begin in the second section by discussing our formal modeling of a DSL, as well as what composition of DSLs within the syntax and semantics of the programming language of the meta-DSL will look like in Coq; and in the third section, we discuss all the mathematical necessities for the algebraic framework of our verification.  



\section{A Formal Model of DSLs and DSL Composition}\label{sec2} 

\subsection{Formalizing DSLs}
The syntax and semantics for the meta-DSL in the maine pipeline for LUMOS will be written in OCaml, while the model for verification will be written in Coq.  Noting Coq can be lifted to OCaml, we expect there to be an interplay between our modeling and the main pipeline for LUMOS.  Our first goal within our verification tasking is to provide a formal model of the syntax and semantics of a DSL in the context of the meta-DSL, since the concept of a DSL is central to patching legacy code.  Within our model of the meta-DSL in Coq, we regard a DSL as a collection of types and a collection of finitary functions (functions of finite arity) on those types from which other finitary functions can be built. The next example illustrates this.  
\begin{example}\label{dslu1}
Let $\text{DSLU} := \stl \stl \text{nat, str}\str,\stl \text{\textbf{fprint}, \textbf{finput}}\str\str$ be the abstraction of a DSL written to patch a function on usernames.  The function $\textbf{fprint}\, n\, str$ takes in a number and a username, and prints the first $n$ letters of the username; the other function, $\textbf{finput}\, str$ takes a username and stores it as a string in memory.  

Within DSLU, we can create a function, denoted by \textbf{firstn}, that prints the first $n$ letters of a function that is written to memory with the composition $\textbf{fprint}\, n \,(\textbf{finput}\, str)$.  
\end{example}

Now to give the Coq formalization of this model, we first specify what the collection of types, $D_T$, in which each DSL, $D$, written by the meta-DSL is.  To do this in Coq, we construct a collection of \textit{type sigils} as an inductive type within Coq, and then provide a function, that acts an embedding, to specify the element of Coq's type system that we assign to each type sigil. For an example, our collection of type sigils could be written as: \verb+ Inductive BNU : Type := b | n | u +.  Then we could write a function in Coq that assigns \verb|b|, \verb|n|, and \verb|u| to the Coq types bool, nat, and unit, respectively. This provides a way to reason about the types when it comes time to evaluate expressions with them.  These two pieces make up what we call our \textit{type universe} (this is what is known as a type universe \`{a} la Tarski, see \cite{M} for details).  In Coq, we use a record type to denote an arbitrary type universe, that we denote by \verb+Univ+.  This is the first component of our definition of a DSL.  

To denote our collection of functions, $D_F$, in our DSL $D$, we specify functions in $D_F$ by their names, by their \textit{type signatures}, and also include the action of the function on \textit{terms in Coq}.  To specify function names, the process is similar to creating a type universe: record the names in Coq as an inductive type (so as \textit{descriptive sigils}), say \verb+funNames : Type+, and then provide an appropriate function in Coq matching sigils in \verb+funNames+ to type signatures.  

The process to specifying a collection of type signatures is also similar to specifying a type universe:  we record the domain of our function as a list of type sigils coming from our type universe, and the range as a type sigil from our type universe.  Our collection of type signatures in Coq is is specified as a record type, denoted by \verb+Sig+.  As a quick example, equality of natural numbers within the type universe \verb+BNU+ would have type signature \verb|n -> n-> b|, and its type signature, \verb|eqNatSig : Sig BNU| would be recorded as two components:  the domain list \verb|[n, n]|, and the return type \verb|b|.    

And lastly, we would have to define the action of our functions on \textit{terms in Coq}, which we do by sending function names to actual function types, using their signatures as input so the signature of the actual function type matches by construction. In Coq, this looks like a function sending a name in \verb|funNames| to a Coq function type created from the signature of the name.

These are the components of our specification of a DSL as record in Coq:  the type universe, as specified by something of type \verb+Univ+, which is \verb+BNU+ in our example; a collection of function names, which are defined by an inductive type in Coq, and would be \verb|funNames| in our example; a collection of type signatures, along with a map that assigns a function name to a given type signature; and the appropriate pairing for function names, type signatures, and term definitions.  All of these components are compiled into a record type in Coq.  In our example, we would call this collection \verb|DSLBNU|, and say that \verb|DSLBNU:DSL| in Coq.  

\subsection{A Model for DSL Composition}
One of the main purposes of the meta-DSL will be to combine abstractions of DSLs that have been written (within the meta-DSL) after processing from preceding modules in the LUMOS pipeline.  Our next step in this part of the verification procedure will be to formalize the notion of combining two DSLs using our formal models for DSLs. Our notion of combining DSLs will be referred to as \textit{gluing DSLs}.

Let $D', D''$ be two DSLs.  A map $\Phi:  D'\ra D''$ consists of a function $\vp:  D'_T\ra D''_T$, such that given $f\in D'_F$, there is an induced function $f_\vp\in D''_F$ that respects the action of $f$ on $D'_T$.  If two DSL abstractions $D'$ and $D''$ are to be composed, and in their composition there are several types that would represent equivalent types (for example, each has a copy of the natural numbers), we want to collapse those copies into a single type in their composition.  To do that, we let $Z$ be a DSL that draws its types from a universal, fixed set $\mct$ (in the LUMOS pipeline this would be the set of all types in OCaml), and first construct maps of DSLs:  $D' \leftarrow Z\ra D''$.  Then the \textit{composition of $D'$ and $D''$ with respect to $Z$}, denoted $D$, is defined by the following commuative square (also known as a \textit{pushout diagram)}:  
\[\begin{tikzcd}
Z \arrow[rr, "\Phi''"] \arrow[dd, "\Phi'"']     &  & D'' \arrow[dd, "h''"]  &   \\
                                       &  &                                         &   \\
D' \arrow[rr, "h'"] &  & D    \\
                                       &  &                                         &
\end{tikzcd}\]
Where the $D, h, h''$ are \textit{universal} with respect to this property, and also make $D$ unique up to unique isomorphism.  

The result of this construction is such that $D_T := D'_T \cup D''_T / \sim_T$, where the equivalence relation $\sim_T$ is induced by the maps:  $D'_T \leftarrow Z_T\ra D''_T$. Similarly for function sets $D'_F$ and $D''_F$, we form $D_F := D'_F \cup D''_F / \sim_F$ by identifying function sets according to the maps $D'_F \leftarrow Z_F\ra D''_F$. 

While this is a general categorical construction, if we want to make a meaningful composition of DSLs within the context of a programming language, we note that $Z$ cannot be taken to be \textit{any} DSL, as demonstrated in the next example.  
\begin{example} \label{exfail}
Consider DSLU from Example \ref{dslu1}, and the DSL, denoted by DSLP, given by $\{\{\verb|int|, \verb|struct|\}, \{\textbf{fprint}, \textbf{ffields}\}\}$. Here the types in $\text{DSLP}_T$ include an integer type \verb|int|, but also a generic \verb|struct| type, which for our purposes will act like records in Coq. The function \textbf{fprint} in DSLP behaves in the same way that \textbf{fprint} in DSLU does, but the function \textbf{ffields} returns the fields of a \verb|struct|. If we try to glue the types \verb|struct| from DSLP and \verb|str| from DSLU together when forming the gluing of DSLU and DSLP, the functionality (or lack thereof) of \textbf{ffields} in the composition may no longer be as intended.  For example, assume we name the new glued type \verb|structstr|, then a problem will arise if there are instances of \verb|structstr| types with \verb|str| inputs and we use \textbf{ffields} on these instances.  
\end{example}

This example illustrates the need for controls on what the DSL $Z$ can be in order to disallow situations in a given context where we may be gluing together types or functions that should not be be glued together.  We provide a a general outline to our approach to this problem:  First, we require a collection of base types, denoted by $\mct$, (which, as noted before, we can regard as the set of types in OCaml), and if two types from the DSLs $D'$ and $D''$, say $t'$ and $t''$, are to be glued together via the DSL $Z$, we require, in the form of proof obligations within Coq, that there is a $t\in \mct$, along with equivalences between $t$ and $t'$, and $t$ and $t''$.  These equivalences must also respect operations within the DSLs $D'$ and $D''$. 

We note that this gives a mathematical description the case in which we have glued together two DSLs along another DSL $Z$.  However, how are we to glue if there are multiple DSLs with more complicated relations?  We discuss the mathematical description of this situation in the next section.  

\section{An Algebraic Framework for our Formal Model} \label{sec3}

\subsection{An Algebraic Model for DSLs}
When we view a DSL with our mathematical definition from Section \ref{sec2} as a collection of a set of types and a set of finitary functions on those types, and add in the ability to compose, in a very unrestricted manner, new functions from old, we find that an algebraic structure that would model this behavior very precisely is a \textit{symmetric colored operad}, which we will just call an \textit{operad}.  The formal definition of an operad is lengthy, so we just include all details that are essential to our paper, and refer to \cite{Y} for the rest.  
\begin{definition}\label{defoperad} 

An \textit{operad} $\mco$ consists of a collection of types, which we will denote by $T$, and for each $n\geq 1$, $d, \und c := c_0,\ldots, c_{n-1}$, a collection of terms $\mco{d\choose \und c}$ and, 
\begin{enumerate}
\item[$\bullet$] for each $c\in T$, an element $\mbbo_c\in\mco{c\choose c}$ called the \textit{$c$-colored unit}; 
\item[$\bullet$] for each $i\in \stl 0,\ldots, n-1\str$, an operator
$$\circ_i: \mco{d\choose \underline c}\times \mco{c_i\choose\underline b} \rightarrow  
\mco{d\choose c_0,\ldots, c_{i-1}, \underline b, c_{i+1},\ldots, c_{n-1}}$$
and axiomatic constraints for associativity, unitary, and symmetry conditions.  

A morphism of operads, $F: \mco\ra \mco'$ consists of a map between types and on terms that commutes with colored units, the operators $\circ_i$, and all axiomatic constraints.  This turns operads into a  category that we denote by $\mbo$.  
\end{enumerate} 

\end{definition} 
While this seems like an extraordinarily abstract definition, the next example helps clarify the roots of the abstraction of Definition \ref{defoperad}.  

\begin{example}\label{exsets}
If we let $T$ be a collection of sets, we can define an operad $\sets_T$ by setting 
$$\sets_T{d\choose c_0,\ldots, c_{n-1}} := \hm(c_0\times\cdots\times c_{n-1}, d)$$
Where the hom-set on the right is the collection of all functions from the set $c_0\times\cdots\times c_{n-1}$ to the set $d$.  In this setting, we can explicitly define the operator $\circ_i$ from Definition \ref{defoperad} which returns, given $f\in \hm(c_0\times\cdots\times c_{n-1}, d) $ and $g\in \hm(b_0\times\cdots\times b_{m-1}, c_i)$, the function $f\circ_i g$ acts on the $(n+m-1)$-tuple $(x_0,\ldots, x_{i-1},\underline y,x_{i+1},\ldots, x_{n-1})$
as:  
$$(f\circ_i g)(x_0,\ldots, x_{i-1},\underline y,x_{i+1},\ldots, x_{n-1}) = f(x_0,\ldots, x_{i-1}, g(\underline y), x_{i+1},\ldots, x_{n-1})$$ 
\end{example}

Now within the formal model of the meta-DSL, it becomes clear that each DSL should be modeled in a similar setting to Example \ref{exsets}.  In this manner, we say the DSL is modeled as an \textit{operad of sets}.  To clarify what we mean, we give an example in a more concrete of setting.  

\begin{example}\label{dslu2}
Let $\text{DSLU} := \stl \stl \text{nat, str}\str,\stl \text{\textbf{fprint}, \textbf{finput}}\str\str$ be as in Example \ref{dslu1}.  We can regular DSLU as an operad of sets, and and use notation from Example \ref{exsets} to say: $\textbf{fprint}\in\text{DSLU}{\text{str}\choose \text{nat, str}}$, and $\textbf{finput}\in\text{DSLU}{\text{str}\choose \text{str}}$.  In particular, we  can write $\textbf{firstn} = \textbf{fprint}\circ_1 \textbf{finput}\in\text{DSLU}{\text{str}\choose \text{nat, str}}$.  
\end{example}

\subsection{An Algebraic Model for the meta-DSL}
At this point, we want each DSL to be represented within the formal algebraic structure of our meta-DSL as an operad of sets.  This leads us to the next question:  What is the proper algebraic representation of the meta-DSL?  As we have noted, one of the main purposes of the meta-DSL will be to compose abstractions of DSLs, so within this framework, we are looking for a mathematical object that provides a way to combine two operads and produces another operad.  

Our approach to this is to use the idea of a \textit{categories of diagrams}.  To start, we give a mathematical definition of a \textit{diagram in a category $\mcc$}.  To wit, if $\mcj$ is a category, a \textit{diagram of shape $\mcj$ in a category $\mcc$} is a covariant functor $D:  \mcj \ra \mcc$.  Moreover, fixing $\mcj$, these objects form a category that we call the \textit{category of diagrams of shape $\mcj$ in $\mcc$}.  We also denote this category by $\mcc^\mcj$, which is standard notation for the category of covariant functors from $\mcj$ to $\mcc$.  The next example is important in our context.  
\begin{example}\label{exspan}
Let $\mcj$ be the category with objects $\bf{-1}, \bf0, \bf1$, and whose non-identity morphisms are given by the diagram $\bf{-1} \leftarrow \bf0 \rightarrow \bf1$.  If $\mcc$ is a category, the image of a diagram $D: \mcj\ra \mcc$ is given by 
$C' \overset{f}{\leftarrow} Z \overset{g}{\rightarrow} C''$ in $\mcc$.  Such a diagram is called a \textit{span}.  
\end{example}

And this is where the mathematics and formalization meet: we can take the \textit{categorical pushout} of a span in Example \ref{exspan} with respect to $Z$.  This has the same properties as the pushout for DSL composition that we described in Section \ref{sec2}, in particular it is unique up to unique isomorphism (see \cite{ML} Ch. III Sec. 3).   

To expand on the idea of combining multiple DSLs along with more complicated relations, we let $\text{Diag}(\mbo)$ be the category in which an object can be a diagram of shape $\mcj$ in $\mbo$ provided $\mcj$ has \textit{finitely many objects}.  The notion of what we are looking for meta-DSL to produce is the \textit{colimit} (see \cite{ML} Ch. III Sec. 3) of an object in $\text{Diag}(\mbo)$.  That is, we can conclude that the meta-DSL can be formalized as the categorical object $\text{Colim}:  \text{Diag}(\mbo)\ra \mbo$.  









%
%
%

\end{document}